\def\draft{n}
\documentclass{amsart}
\usepackage[headings]{fullpage}
\usepackage{amssymb,epic,eepic,epsfig}
\usepackage{graphicx}

\theoremstyle{plain}

\newtheorem{theorem}{Theorem}

\newtheorem{proposition}{Proposition}[section]
\newtheorem{lemma}[proposition]{Lemma}
\newtheorem{corollary}[proposition]{Corollary}

\newtheorem{conjecture}{Conjecture}

\theoremstyle{definition}
\newtheorem{definition}[proposition]{Definition}
\newtheorem{problem}{Problem}
\newtheorem{question}{Question}

\theoremstyle{remark}
\newtheorem{example}[proposition]{Example}

\newtheorem{remark}[proposition]{Remark}

\def\printname#1{
        \if\draft y
                \smash{\makebox[0pt]{\hspace{-0.5in}
                        \raisebox{8pt}{\tt\tiny #1}}}
        \fi
}

\newcommand{\psdraw}[2]
         {\begin{array}{c} \hspace{-1.3mm}
        \raisebox{-4pt}{\epsfig{figure=draws/#1.eps,width=#2}}
        \hspace{-1.9mm}\end{array}}

\newlength{\standardunitlength}
\catcode`\@=11
\long\def\@makecaption#1#2{%
     \vskip 10pt

\setbox\@tempboxa\hbox{
       \small\sf{\bfcaptionfont #1. }\ignorespaces #2}%
     \ifdim \wd\@tempboxa >\captionwidth {%
         \rightskip=\@captionmargin\leftskip=\@captionmargin
         \unhbox\@tempboxa\par}%
       \else
         \hbox to\hsize{\hfil\box\@tempboxa\hfil}%
     \fi}
\font\bfcaptionfont=cmssbx10 scaled \magstephalf
\newdimen\@captionmargin\@captionmargin=2\parindent
\newdimen\captionwidth\captionwidth=\hsize
\catcode`\@=12

\def\lbl#1{\label{#1}\printname{#1}}

\def\BN{\mathbb N}
\def\BZ{\mathbb Z}

\def\BQ{\mathbb Q}
\def\BR{\mathbb R}
\def\BC{\mathbb C}

\def\D{\Delta}

\def\a{\alpha}

\def\s{\sigma}

\def\ga{\gamma}

\def\la{\langle}
\def\ra{\rangle}

\def\e{\epsilon}

\def\De{\Delta}
\def\d{\delta}
\def\b{\beta}

\def\s{\sigma}

\def\longto{\longrightarrow}

\def\SL{\mathrm{SL}}
\def\pt{\partial}

\def\ep{\epsilon}
\def\deg{\mathrm{deg}}
\def\js{\mathrm{js}}
\def\bs{\mathrm{bs}}
\def\jd{\mathrm{jd}}
\def\qs{\mathrm{qs}}

\begin{document}


\title[The Jones slopes of a knot]{
The Jones slopes of a knot}
\author{Stavros Garoufalidis}
\address{School of Mathematics \\
         Georgia Institute of Technology \\
         Atlanta, GA 30332-0160, USA \\ 
         {\tt http://www.math.gatech} \newline {\tt .edu/$\sim$stavros } }

\thanks{The author was supported in part by NSF. \\
\newline
1991 {\em Mathematics Classification.} Primary 57N10. Secondary 57M25.
\newline
{\em Key words and phrases: Knot, link, Jones polynomial, Jones slope,
Jones period, quasi-polynomial, alternating knots, signature, pretzel knots,
polytopes, Newton polygon, incompressible surfaces, slope, Slope Conjecture.
}
}

\date{May 21, 2010} 


\begin{abstract}
The paper introduces the Slope Conjecture which relates the
degree of the Jones polynomial of a knot and its parallels with the slopes
of incompressible surfaces in the knot complement. More precisely, we 
introduce two knot invariants, the Jones slopes (a finite set of rational 
numbers) and the Jones period (a natural number) of a knot in 3-space. 
We formulate a number of conjectures for these invariants and verify them
by explicit computations for the class of alternating knots, the knots
with at most 9 crossings, the torus knots and the $(-2,3,n)$ pretzel knots.
\end{abstract}

\maketitle

\tableofcontents

\section{Introduction}
\lbl{sec.intro}

\subsection{The degree of the Jones polynomial and incompressible
surfaces}
\lbl{sub.cj}

The paper introduces an explicit conjecture relating the 
degree of the Jones polynomial of a knot (and its parallels) with
slopes of incompressible surfaces in the knot complement. We give an
elementary proof of our conjecture for alternating knots and torus knots, 
and check it with explicit computations for non-alternating knots with 
$8$ and $9$ crossings, and for the $(-2,3,p)$ pretzel knots.

One side of our conjecture involves the growth rate of the degree 
$\d_K(n)$ (with respect to $q$) of the colored Jones function $J_{K,n}(q)
\in \BZ[q^{\pm 1}]$ 
of a knot. The other side involves the finite set $\bs_K$ of {\em slopes 
of incompressible, $\pt$-incompressible orientable
surfaces} in the complement of $K$,
where the slopes are normalized so that the longitude has slope $0$
and the meridian has slope $\infty$; \cite{Ha}. 
To formulate our conjecture, we need a definition. Recall that $x \in \BR$
is a cluster point of a sequence $(x_n)$ of real numbers if for every $\ep>0$
there are infinitely many indices $n \in \BN$ such that $|x-x_n| < \ep$. 
Let $\{x_n\}'$ denote the set of {\em cluster points} of a sequence $(x_n)$.

\begin{definition}
\lbl{def.js}
\rm{(a)}
For a knot $K$, define the {\em Jones slopes} $\js_K$ by:
\begin{equation}
\lbl{eq.js}
\js_K=\{ \frac{2}{n^2}\deg(J_{K,n}(q)) \,\, | n \in \BN \}'
\end{equation}
\rm{(b)} Let $\bs_K$ denote the set of boundary slopes of incompressible 
surfaces of $K$.
\end{definition}
A priori, the structure and the cardinality of the set $\js_K$ is not
obvious. On the other hand, 
it is known that $\bs_K$ is a finite subset of $\BQ \cup \{\infty\}$;
see \cite{Ha}. Normal surfaces are of special interest because of their
relation with {\em exceptional} Dehn surgery, and the $\mathrm{SL}(2,\BC)$
character variety and hyperbolic geometry, see for
example \cite{Bu,CGLS,CCGLS,KR,LTi,Mv}. 

\begin{conjecture}
\lbl{conj.1}(The Slope Conjecture)
For every knot we have
\begin{equation}
2 \js_K \subset \bs_K.
\end{equation}
\end{conjecture}
Before we proceed further, and to get a better intuition about this
conjecture, let us give three illustrative examples.

\begin{example}
\lbl{ex.817}
For the alternating knot $8_{17}$ we have:
\begin{eqnarray*}
\d(n) &=& 2n^2+2n
\\
\d^*(n) &=& -2n^2-2n
\end{eqnarray*}
where $\d_K(n)$ and $\d^*_K(n)$ are the maximum and the the minimum degree of 
$J_{K,n}(q)$ with respect to $q$.
On the other hand, according to \cite{Cu}, the Newton polygon 
(based on the geometric component of the character variety) 
has $44$ sides and its slopes (excluding multiplicities) are
$$
\{-14,-8,-6,-4,-2,0,2,4,6,8,14,\infty\}.
$$
The reader may observe that $\d(n)$ and $\d^*(n)$ are quadratic polynomials
in $n$ and four times the leading terms of
$\d(n)$ and $\d^*(n)$ are boundary slopes (namely $8$ and $-8$), and moreover 
they agree with $2 c^+$ and $-2 c^-$ where $c^{\pm}$ is the number of 
positive/negative crossings of $8_{17}$. In addition, as Y. Kabaya
observed (see \cite{Ka}), $8_{17}$ has slopes $\pm 14$
outside the interval $[-2c^-,2c^+]=[-8,8]$.
\end{example}

\begin{example}
\lbl{ex.237}
For the non-alternating pretzel knot $(-2,3,7)$ we have:
\begin{eqnarray*}
\d(n) &=& \left[\frac{37}{8} n^2 + \frac{17}{2} n \right]=
\frac{37}{8} n^2 + \frac{17}{2} n + \e(n),
\\
\d^*(n) &=& 5n
\end{eqnarray*}
where $\e(n)$ is a periodic sequence of period $4$ given by $0,1/8,1/2,1/8$
if $n\equiv 0,1,2,3 \bmod 4$ respectively. $(-2,3,7)$ is a Montesinos knot
and its boundary slopes are given by
$$
\{0,16, 37/2, 20\}
$$
(see \cite{HO} and \cite{Du} and compare also with \cite{Ma}). 
In this case, $\d(n)$ is no longer a quadratic polynomial of $n$. Instead 
$\d(n)$ is a quadratic quasi-polynomial with fixed leading term $37/8$. 
Moreover, four times this leading term is a slope of the knot. This number
was the motivating example that eventually lead to the results of this paper.
Likewise, four times 
the $n^2$-coefficient of $\d^*(n)$ is $0$, which is also a boundary slope 
of this knot.
\end{example}

\begin{example}
\lbl{ex.mutation}
The pretzel knots $(2,3,5,5)$ and $(2,5,3,5)$ are mutant, alternating
and Montesinos. Since they are mutant, their colored Jones functions 
(thus their Jones slopes) agree; see \cite{Kf}. Since they are alternating,
their common Jones slopes are $\js=c^+=15$ and $\js^*=-c^-=0$; see Theorem
\ref{thm.1alt}. Since they are Montesinos, Hatcher-Oertel's algorithm
implemented by Dunfield (see \cite{Du}) implies that their boundary
slopes are given by
\begin{eqnarray*}
\bs_{(2,3,5,5)}&=&
\{0, 4, 6, 10, 14, 16, 18, 20, 22, 24, 26, 28, 30, \infty\} \\
\bs_{(2,5,3,5)}&=& \{0, 4, 6, 10, 14, 16, 18, 20, 22, 26, 28, 30, \infty\}
\end{eqnarray*}
It follows that the set $\bs_K$ is not invariant under knot mutation. 
\end{example}

\subsection{The degree of the colored Jones function is a quadratic
quasi-polynomial}
\lbl{sub.degree}

In the previous section, we took the shortest path to formulate a 
conjecture relating the degree of the colored Jones function of a knot
with incompressible surfaces in the knot complement. In this section
we will motivate our conjecture, and add some {\em structure} to it.
Let us recall that in 1985, Jones introduced the famous {\em Jones polynomial}
$J_K(q) \in \BZ[q^{\pm 1}]$ of a knot (or link) $K$ in 3-space; \cite{Jo}. 
The Jones polynomial of a knot is a Laurent polynomial with integer 
coefficients that tightly encodes information about the topology and the 
geometry of the knot. 

Unlike the Alexander polynomial, not much is known about the topological
meaning of the coefficients of the Jones polynomial, nor about its degree,
nor about the countable set of Jones polynomials of knots. 

This unstructured behavior of the Jones polynomial becomes more
structured when one fixes a knot $K$ and considers a stronger invariant, 
namely the {\em colored Jones function} $J_{K,n}(q)$. The latter is a 
sequence of elements of $\BZ[q^{\pm 1}]$ indexed by $n \in \BN$ which encodes 
the TQFT invariants of a knot colored by the irreducible $(n+1)$-dimensional 
representation of $\mathrm{SU}(2)$, and normalized to be $1$ at the unknot; 
see \cite{Tu2}. With these conventions, $J_{K,0}(q)=1$ and $J_{K,1}(q)$ is the 
Jones polynomial of $K$.

In many ways, the sequence $J_{K,n}(q)$ is better behaved, and
suitable limits of the sequence $J_{K,n}(q)$ have a clear 
topological or geometric meaning. Let us give three instances of this
phenomenon:

\begin{itemize}
\item[(a)]
A suitable formal power series limit $J_{K,n-1}(e^h) \in \BQ[[h,n]]$ (known 
as the Melvin-Morton-Rozansky Conjecture) equals to $1/\D_K(e^{n h})$
and determines the {\em Alexander polynomial} 
$\D_K(t)$ of $K$ (see \cite{B-NG}).
\item[(b)]
An analytic limit $J_{K,N-1}(e^{\a/N})$ for small complex numbers $\a$ near
zero equals to $1/\D_K(e^{\a})$ 
also determines the Alexander polynomial of $K$; see 
\cite[Thm.2]{GL1}.
\item[(c)]
The exponential growth rate of the sequence 
$J_{K,N-1}(e^{2 \pi i/N})$ (the so-called {\em Kashaev invariant})
of a hyperbolic knot is conjectured to equal to the {\em volume} of $K$,
divided by $2 \pi$; see \cite{Ks}.
\end{itemize}
On the other hand, one can easily construct hyperbolic knots with equal 
Jones polynomial but different Alexander polynomial and volume. 

Some already observed structure regarding the colored Jones function
$J_{K,n}(q)$ is that it is $q$-holonomic, i.e., it satisfies a linear
recursion relation with coefficients in $\BZ[q^n,q]$; see \cite{GL1}.
The present paper is concerned with another notion of regularity, namely
the degree of the colored Jones function $J_{K,n}(q)$ with respect to $n$.
Since little is known about the degree of the Jones polynomial of a knot,
one might expect that there is little to say about the degree of the colored 
Jones function $J_{K,n}(q)$. Once observed, the regularity of the
degree seems obvious as Bar-Natan suggests; see \cite[Lem.3.6]{B-NL}
and \cite{Me}. Moreover, the degree of the colored Jones function
motivates the introduction of two knot invariants, the Jones
slopes of a knot (a finite set of rational numbers) and the Jones period
of a knot (a natural number).

\subsection{$q$-holonomic functions and quadratic quasi-polynomials}
\lbl{sub.qholo}
 
To formulate our new notion of regularity, we need to recall
what is a quasi-polynomial. A {\em quasi-polynomial} $p(n)$ is a 
function 
$$p: \BN \longto \BN, \qquad p(n)=\sum_{j=0}^d c_j(n) n^j
$$ 
for some $d \in \BN$ where $c_j(n)$ is a periodic fuction with integral 
period for $j=1,\dots,d$; \cite{St,BR}.
If $c_d(n)$ is not identically zero, then the degree of $p$ is $d$.
We will focus on two numerical invariants of a quasi-polynomial, its
period and its slopes.  

\begin{definition}
\lbl{def.quasi}
\rm{(a)} The {\em period} $\pi$ of a quasi-polynomial
$p(n)$ as above is the common period of $c_j(n)$. 
\newline
\rm{(b)} The set of {\em slopes} of a quadratic quasi-polynomial 
$p(n)$ is the finite set of twice the
rational values of the periodic function $c_2(n)$.
\end{definition}
Notice that if $p(n)$ is a quasi-polynomial of period $\pi$, then there 
exist polynomials $p_0, \dots, p_{s-1}$ such that $p(n) = p_i(n)$ when 
$n \equiv i \bmod \pi$, and vice-versa. Notice also that the set of slopes of
a quasi-polynomial is always a finite subset of $\BQ$.

Quasi-polynomials of period $1$ are simply polynomials. Quasi-polynomials
appear naturally in counting problems of lattice points in rational
convex polytopes; see for example \cite{BP,BR,BV,Eh,St}. 
In fact, if $P$ is a rational convex polytope, then the number of lattice 
points of $nP$ is the so-called {\em Ehrhart} quasi-polynomial of $P$, 
useful in many enumerative questions \cite{BP,BR,BV,Eh,St}.

The next theorem seems obvious, once observed. The proof, given in
\cite{Ga2}, uses ideas from {\em differential Galois theory} of $D$-modules 
and the key {\em Skolem-Mahler-Lech theorem} from number theory. 
Let $\deg(f(q))$ denote the degree of a rational function $f(q)$ 
with respect to $q$.

\begin{theorem}
\lbl{thm.Ga}\cite{Ga2}
If $f_n(q)$ is a $q$-holonomic sequence of rational functions, then for 
large $n$, $\deg(f_n(q))$ is a quadratic quasi-polynomial. Moreover,
the leading term of $f_n(q)$ satisfies a linear recursion relation with
constant coefficients.
\end{theorem}

The restriction {\em for large $n$} in Theorem \ref{thm.Ga} is necessary,
since the sequence $((1+(-1)^n)q^{n^2}+q^{17})$ is $q$-holonomic, and its degree 
(given by $17$ if $n \leq 4$, by $n^2$ if $n \geq 5$ is even and
by $17$ if $n \geq 5$ is odd) is not a quasi-polynomial. 
On the other hand, if $n \geq 5$, the degree is given by the quadratic
quasi-polynomial $n^2 (1+(-1)^n)/2 + 17$.

\begin{corollary}
\lbl{cor.jslopes}
If $f_n(q)$ is a $q$-holonomic sequence of rational functions
and $\d(n)=\deg(f_n(q))$ then
$$
\{ \frac{2}{n^2}\d(n) \,\, | n \in \BN \}'=\mathrm{slopes}(\d)
$$
\end{corollary}

\begin{proof}
Consider a subsequence $k_n$ of the natural numbers such that
$$
l=\lim_n 2 \frac{\d(k_n)}{k_n^2}
$$
Since $d(n)$ is (for large $n$) given by a quasi-polynomial, it follows that
$d(n)=c_2(n) n^2 + O(n)$ for a periodic function $c_2: \BN \longto \BQ$
and for all $n \in \BN$.
Since $c_2$ takes finitely many values, it follows that there is a subsequence
$m_n$ of $k_n$ such that $c_2(m_n)=s$ for all $n$, where $s$ is a slope of 
$\d(n)$. Since $d(n)=c_2(n) n^2 + O(n)$ for all $n$, it follows that
$l=2 s$. I.e., $l$ is a slope of $\d$. Conversely, it is easy to see that
every slope of $\d$ is the limit point for a subsequence taken to be
an arithmetic progression on which $c_2$ takes a constant value.
\end{proof}

\subsection{The Jones slopes and the Jones period of a knot}
\lbl{sub.jslopes}

Given a knot $K$, we set
\begin{equation}
\d_K(n)=\deg(J_{K,n}(q))
\end{equation}
Combining the $q$-holonomicity of the colored Jones function $J_{K,n}(q)$
of a knot $K$ with Theorem \ref{thm.Ga}, it follows that $\d_K$ is a quadratic
quasi-polynomial for large $n$.

\begin{definition}
\lbl{def.jslopes}
\rm{(a)} The {\em Jones period} $\pi_K$ is the period of $\d_K$.
\newline
\rm{(b)} The {\em Jones slopes} $\js_K$ is the finite set of slopes of 
$\d_K$.
\end{definition}
Corollary \ref{cor.jslopes} implies that
for every knot $K$ we have:
$$
\{ \frac{2}{n^2}\deg(J_{K,n}(q)) \,\, | n \in \BN \}'=\mathrm{slopes}(\d_K)
$$
where the right-hand side (and consequently, the left-hand side, too)
is a finite subset of $\BQ$.

\begin{lemma}
\lbl{lem.jslopes}
For every knot $K$ we have:
$$
\{ \frac{2}{n^2}\deg(J_{K,n}(q)) \,\, | n \in \BN \}'=\mathrm{slopes}(\d_K)
$$
where $\mathrm{slopes}(\d_K)$ is a finite subset of $\BQ$.
\end{lemma}

The $\d_K$ invariant records the growth rates of the (maximum) degree of the 
colored Jones function of $K$. We can also record the minimum degree as 
follows. If $K^*$ denotes the mirror image of $K$, then 
$J_{K^*,n}(q)=J_K(q^{-1})$. Let us define
\begin{equation}
\lbl{eq.Jstar}
\d^*_K(n)=-\d_{K^*}(n)=\text{mindeg}(J_{K,n}(q)), \qquad
\js^*_K=\text{slopes}(\d^*_K).
\end{equation}
Notice that Conjecture \ref{conj.1} applied to $K^*$ implies that
$$
\js^*_K \subset \bs_K.
$$

The next proposition gives a bound for the Jones slopes of a knot $K$ in terms
of the number $c^{\pm}_K$ of positive/negative crossings of a planar 
projection.

\begin{proposition}
\lbl{prop.c+-}
For every knot $K$, every $s \in \js_K$ and every $s^* \in \js^*_K$ we have
\begin{equation}
\lbl{ss*}
-c^-_K \leq s^*,  s \leq c^+_K.
\end{equation}
\end{proposition}
The reader may compare the above lemma with Example \ref{ex.817}.

Our next theorem confirms Conjecture \ref{conj.1} for all alternating knots.
Consider a {\em reduced planar projection} of $K$ with $c^{\pm}_K$ crossings 
of positive/negative sign.

\begin{theorem}
\lbl{thm.1alt}
If $K$ is alternating, then
\begin{equation}
\lbl{eq.alt1}
\pi_K=1, \qquad \js_K=\{c^+_K\}, \qquad \js^*_K=\{-c^-_K\}.
\end{equation}
In addition, the two checkerboard surfaces of $K$ (doubled, if need, to
make them orientable) are incompressible
with slopes $2c^+_K$ and $-2c^-_K$.
\end{theorem}

Our final lemma relates the Jones slopes and the period of a knot.

\begin{lemma}
\lbl{lem.spi}
If $a(n)$ is an integer-valued quadratic quasi-polynomial with period $\pi$, 
then for every slope $s$ of $a(n)$ we have
$$
s \pi^2 \in \BZ.
$$
In particular, for every knot $K$ we have
$$
\pi_K^2 \,\, \js_K \subset \BZ, \qquad \pi_K^2 \,\, \js^*_K \subset \BZ.
$$
Thus, if a knot has a non-integral Jones slope, then it has period bigger 
than $1$.
\end{lemma}

\subsection{The symmetrized Jones slopes and the signature of a knot}
\lbl{sub.alt}

In this section we discuss the symmetrized version $\d^{\pm}_K$ of $\d_K$:
\begin{equation}
\lbl{eq.delta}
\d^+_K=\d_K-\d^*_{K}, \qquad \d^-_K=\d_K+\d^*_{K}.
\end{equation}
Of course, $\d_K=1/2(\d^+_K+\d^-_K)$ and $\d^*_K=1/2(-\d^+_K+\d^-_K)$. Pictorially,
we have:
$$
\psdraw{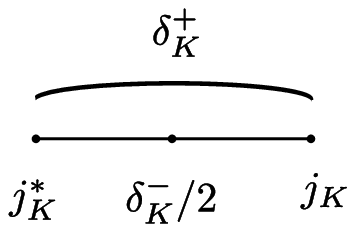}{1.3in}
$$
As we will see below, the symmetrized degree $\d^{\pm}_K$ of the colored
Jones function has a different flavor, and relates (at
least for alternating knots) to the signature of the knot.
$\d^+_K(n)$ is the {\em span} of $J_{K,n}(q)$, i.e., the difference between
the maximum and minimum degree of $J_{K,n}(q)$. On the other hand, $\d^-K(n)$ 
is the sum of the minimum and maximum degree of $J_{K,n}(q)$, and appears
to be less studied. Of course, $\d^{\pm}_K$ are quadratic quasi-polynomials.
Since the colored Jones function is multiplicative under 
connected sum, and reverses $q$ to $q^{-1}$ under mirror image, it follows that

\begin{equation}
\lbl{eq.mirrordelta}
\d^-_{K_1 \# K_2}=\d^-_{K_1}+ \d^-_{K_2}, \qquad \d^-_{K^*}=-\d^-_K.
\end{equation}
Our next theorem computes the $\d^{\pm}_K$ quasi-polynomials 
of an alternating knot $K$ in terms of three basic invariants: 
the {\em signature} $\s_K$, the {\em writhe} $w_K$
and the {\em number of crossings} $c_K$ of a reduced projection of $K$. 
Our result follows from elementary linear algebra using the results of 
Kauffman, Murasugi and Thistlethwaite, \cite{Kf,Mu,Th}, 
further simplified by Turaev \cite{Tu1}. See also \cite[p.42]{Li} and 
\cite[Prop.2.1]{Le}.

\begin{theorem}
\lbl{thm.2alt}
\rm{(a)} For all alternating knots $K$ we have:
\begin{eqnarray}
\lbl{eq.deltap}
\d^-_K(n) &=& \frac{w_K}{2} n^2 + \frac{w_K-2 \s_K}{2} n \\
\lbl{eq.deltan}
\d^+_K(n) &=& \frac{c_K}{2} n^2 + \frac{c_K}{2} n 
\end{eqnarray}
\rm{(b)} The Jones polynomial $J_K(q)$ determines $c_K$ by:
\begin{equation}
\lbl{eq.a3}
\d^+_K(1)=c_K.
\end{equation}
\rm{(c)} The Jones polynomial of $K$ and its 2-parallel  
determines $w_K$ and $\s_K$ by:
\begin{equation}
\lbl{eq.a4}
\s_K=-3\d^-_K(1)+\d^-_K(2), \qquad w_K=-2\d^-_K(1)+\d^-_K(2).
\end{equation}
\end{theorem}

\begin{remark}
\lbl{rem.2}
Part (c) of Theorem \ref{thm.2alt} is sharp.
The Jones polynomial of an alternating knot determines the number of crossings,
but it does not determine the signature nor the writhe of the knot. 
Shumakovitch provided us with a table of pairs of 
alternating knots with up to $14$ crossings (using the Thistlethwaite notation)
with equal Jones polynomials and unequal signature. An example of such a pair 
is with $12$ crossings is the knot $12a_{669}$ and its mirror image:
$$
\psdraw{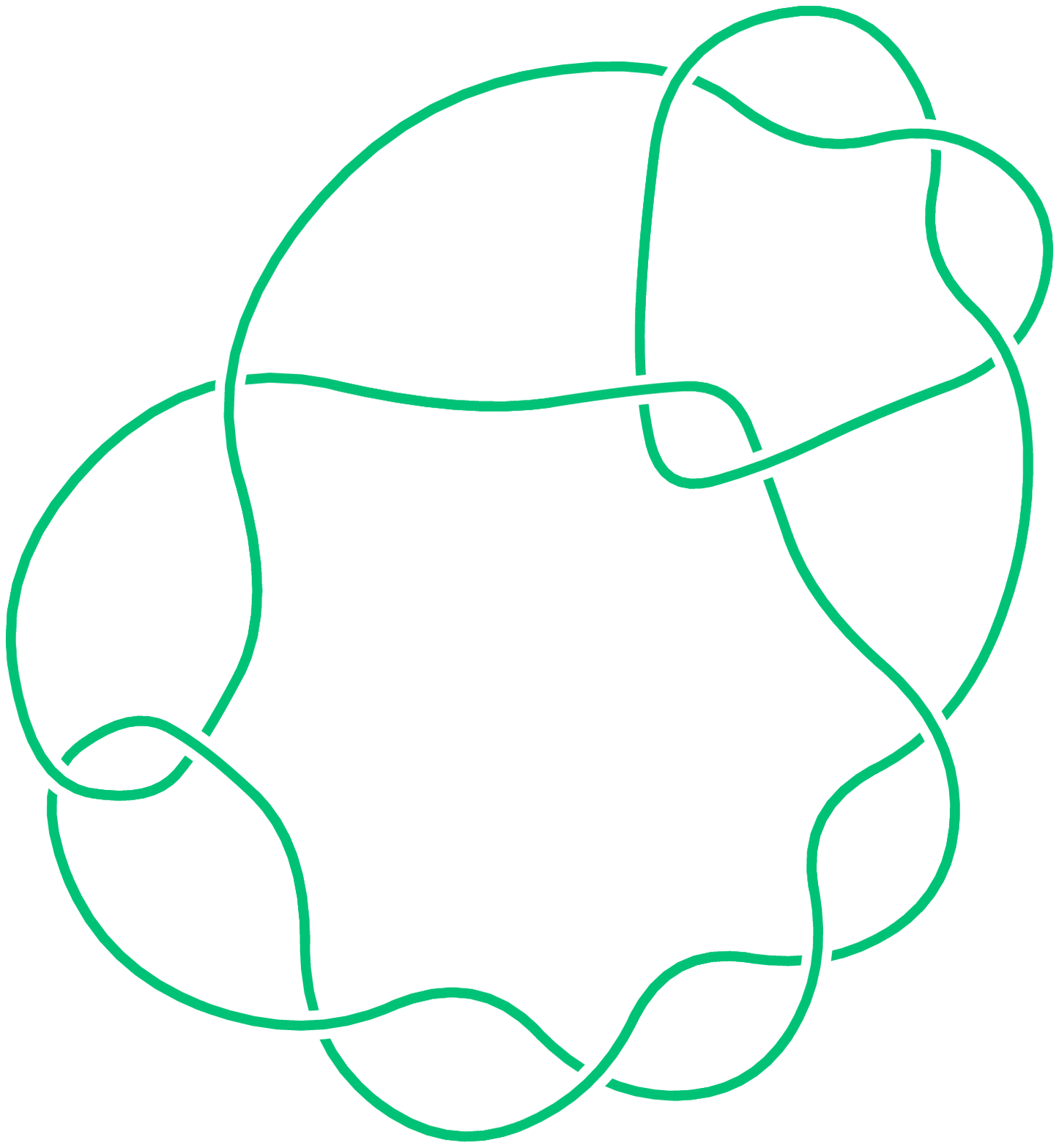}{1.5in}
$$ 
$12a_{669}$ has Jones polynomial
$$
J_{12a_{669},1}(q)
=- \frac{1}{q^6} + \frac{2}{q^5} - \frac{4}{q^4} + \frac{6}{q^3} 
- \frac{7}{q^2} + \frac{9}{q} -9 + 9 q - 7 q^2 + 
6 q^3 - 4 q^4 + 2 q^5 - q^6
$$
and signature $-2$. Since the signature is nonzero, $12a_{669}$ is {\em not
amphicheiral}, and yet has {\em palindromic} Jones polynomial. The next colored
Jones polynomial is given by:
{\small
\begin{eqnarray*}
J_{12a_{669},2}(q)&=&
\frac{1}{q^{17}} - \frac{2}{q^{16}} + \frac{1}{q^{15}} + \frac{3}{q^{14}} 
- \frac{7}{q^{13}} 
+ \frac{4}{q^{12}} + \frac{6}{q^{11}} - \frac{13}{q^{10}} + \frac{6}{q^{9}} 
+ \frac{8}{q^8} 
- \frac{15}{q^7} + \frac{7}{q^6} + \frac{7}{q^5} - \frac{15}{q^4} 
+ \frac{11}{q^3} + 
\frac{6}{q^2} - \frac{21}{q} 
\\ & &
+ 18 + 9 q - 30 q^2 + 20 q^3 + 15 q^4 - 35 q^5 
+ 16 q^6 +  20 q^7 - 32 q^8 + 7 q^9 + 22 q^{10} - 22 q^{11} - 2 q^{12} 
+ 17 q^{13} 
\\ & &
-  9 q^{14} - 5 q^{15} + 7 q^{16} - q^{17} - 2 q^{18} + q^{19}
\end{eqnarray*}  
}
\noindent
and it is far from being palindromic.
This is another example where the pattern of the Jones polynomial is
blurred, but the pattern of the colored Jones function is clearer.
\end{remark}

Results similar to Theorem \ref{thm.1alt} and \ref{thm.2alt} have been
also been obtained independently in \cite{CT} using the Jones polynomial 
of an alternating knot.

Let us end this section with two comments. In this paper
$K$ is a knot, but without additional effort one can state 
similar results for a link (and even a knotted trivalent graph, or
quantum spin network) in 3-space. In addition, 
we should point out that there are deeper 
aspects of stability and integrality of the coefficients of the 
colored Jones function. We will discuss them in a future publication.

\subsection{Plan of the proof}
\lbl{sub.plan}
In Section \ref{sec.alt} we use the Kauffman bracket skein module 
and the work of Kauffman, Le, Murasugi Thistlethwaite to give a proof
of Proposition \ref{prop.c+-} and Theorems \ref{thm.1alt} and 
\ref{thm.2alt}. This proves our Slope Conjecture for all alternating knots.

In Section \ref{sec.nonalt} we give computational evidence for the
degree of the colored Jones polynomials (and for the Slope Conjecture)
of the non-alternating knots with $8$ and $9$ crossings. In addition,
we verify the Slope Conjecture for all $(-2,3,p)$ pretzel knots 
(using the fusion state-sum formulas of the colored Jones function
studied in \cite{GL3,Co,GVa}), and all torus knots (using Morton's
formula for the colored Jones function of those, \cite{Mo}). The reader
may note that the pretzel knots $(-2,3,p)$ are well-known examples of
Montesinos knots, with non-integral slopes (for $p \neq -1,1,3,5$), which
are {\em not} quasi-alternating, thus the results of \cite{FKP} do not
apply for this family.

\section{Future directions}
\lbl{sec.future}

In this section we will discuss some future directions, and pose some
questions, problems and conjectures.

The Slope Conjecture 
involves knots in 3-space, and relates the degree of the colored Jones
function to their set of boundary slopes. The Slope Conjecture may be
extended in three different directions: one may consider (a) links in 3-space,
(b) general 1-cusped manifolds and (c) general Lie algebras. These
extensions have been considered by the author. In \cite{GVu} we introduce
a Slope Conjecture for knots and arbitrary Lie algebras. One may also
consider 1-cusped manifolds with the homology of $S^1 \times D^2$, i.e.,
knot complements in integer homology spheres. For those, the colored Jones
polynomial exists (though it is {\em not} a polynomial, but rather an
element of the Habiro ring; see \cite{GL2}). The colored Jones function is
$q$-holonomic, and we can use the non-commutative $A$-polynomial of this
sequence to define the set of slopes. 
 
The following problem appears mysterious and tandalizing.

\begin{problem}
\lbl{prop.1}
Understand the {\em Selection Principle} which selects 
\begin{itemize}
\item[(a)]the Jones slope
from the set of boundary slopes, 
\item[(b)]
the colored Jones function from the vector space of solutions
of the linear $q$-difference equation.
\end{itemize}
\end{problem}

Now, let us post some questions, based on the limited experimental
evidence. 
Our first question concerns the class of knots of Jones period $1$.

\begin{question}
\lbl{que.1}
Is it true that a prime knot $K$ is alternating if and only if $\pi_K=1$?
\end{question}
\noindent
Theorem \ref{thm.1alt} proves the easy direction of the above conjecture.
As a warm-up for our next question, 
several authors have studied the {\em diameter} $d_K$ of the set $\bs_K$
$$
d_K=\max\{|s-s'| \,\,\, s,s' \in \bs_K \}
$$ 
See for example \cite{IM1,IM2,MMR}. 
Y. Kabaya pointed out to us that there are alternating knots $K$ with
diameter bigger than twice the number of crossings; see Example 
\ref{ex.817}. Let $\jd_K$ denote the {\em Jones diameter}

\begin{equation}
\lbl{eq.jdiam}
\jd_K=\max\{|s-s^*| \,\,\, s \in \js_K, \,\, s^* \in \js^*_K \}.
\end{equation}
Proposition \ref{prop.c+-} shows that $\jd_K \leq c_K$. Moreover, the
bound is achieved for alternating knots, and more generally for adequate
knots; see \cite{FKP}. The next question
concerns the class of knots of maximal Jones diameter.

\begin{question}
\lbl{que.2}
Is it true that a prime knot $K$ is adequate if and only if $\jd_K=c_K$?
\end{question}
The class of alternating knots is included in two
natural classes: {\em quasi-alternating knots}, and {\em adequate knots}.
Knot Homology (and its exact triangles) 
can tell whether a knot is quasi-alternating or not
(see \cite{OM}), but it seems hard to tell whether a knot is alternating or 
not. Adequate knots appeared in \cite{LT} in relation to the Jones polynomial 
and also in \cite{FKP}. It was pointed out to us by D. Futer, E. Kalfagianni
J. Purcell and P. Ozsv\'ath that the pretzel knots $(-2,3,p)$ are not 
adequate nor quasi-alternating for $p>5$.

In all examples, the set $\js_K$ consists of a single
element, whereas the set $\bs_K$ can have arbitrarily many elements. Thus,
Conjecture \ref{conj.1} sees only a small part of the set $\bs_K$.
Our next conjecture claims that the colored Jones function $J_{K,n}(q)$
of $K$ may see all the elements in $\bs_K$. To formulate it,
recall that $J_{K,n}(q)$ is a $q$-holonomic sequence, and satisfies 
a unique, minimal order recursion relation of the form
$$
\sum_{k=0}^d a_k(q^n,q) J_{K,n+k}(q)=0
$$
where $a_k[u,v] \in \BQ[u,v]$ are polynomials with greatest common divisor
$1$;  see \cite{Ga1}. The 3-variable polynomial 
$qA_K(E,Q,q)=\sum_{k=0}^d a_k(Q,q) E^k$ is often called the {\em non-commutative}
$A$-polynomial of $K$. The AJ Conjecture of \cite{Ga1} 
states that every irreducible factor of $qA_K(L,M,1)$ is a factor of 
$A_K(L,M^2)$ or is $L$-free, and vice-versa. Here $A_K$ denotes the 
$A$-polynomial of $K$ and 
$A_K$ contains all components of the $\SL(2,\BC)$ character variety
of $K$ (including the abelian one). Let $\bs^A_K$ denote the slopes of
the Newton polygon of $A$. These are the so-called {\em visible slopes}
of a knot. It follows by Culler-Shalen theory 
(see \cite{CS,CGLS,CCGLS}) that $\bs^A_K \subset \bs_K$.

Let us define the $q$-{\em Newton polytope} $qN_K$ of $K$ to be the convex 
hull of the monomials $q^cQ^bE^a$ of $qA_K(E,Q,q)$. $qN$ is a convex polytope 
in $\BR^3$, and we may consider the image of it in $\BR^2$ under the projection
map $\BR^3 \longto \BR^2$ which maps $(a,b,c)$ to $(a,b)$ (i.e., sends
the monomial $q^cQ^bE^a$ to $Q^bE^a$). 

\begin{definition}
\lbl{def.qslopes}
The $q$-{\em slopes} $\qs_K$ of a knot $K$ are the slopes of the projection
of $qN$ to $\BR^2$.
\end{definition}
 
\begin{problem}
\lbl{prob.2}
Show that for every knot $K$ we have 
\begin{equation}
\lbl{eq.qsbs}
2 \qs_K = \bs^A_K.
\end{equation}
\end{problem}
It is easy to see that for every knot $K$ we have $\js_K \subset \qs_K$. 
In fact, this holds for arbitrary $q$-holonomic sequences; see 
\cite[Prop.1.2]{Ga4} for a detailed discussion. Thus, Problem 
\ref{prob.2} implies Conjecture \ref{conj.1}. The AJ Conjecture motivates
Problem \ref{prob.2}. This is discussed in detail in \cite{Ga4}.

Our next problem concerns the symmetrized quasi-polynomial $\d^-$ of a knot 
from \eqref{eq.delta}. Although $\d^-$ is not a concordance invariant, it
determines the signature of an alternating knot.

\begin{problem}
\lbl{prob.3}
Show that $\d^-$ determines a Knot Homology invariant.
\end{problem}

\section{The Jones slopes and the Jones period of an alternating knot}
\lbl{sec.alt}

In this Section we prove Proposition \ref{prop.c+-}
and Theorems \ref{thm.1alt} and \ref{thm.2alt}
for an alternating knot $K$, using the Kauffman bracket presentation
of the colored Jones polynomial.

The following lemma of Le \cite[Prop.2.1]{Le},
(based on well-known properties of the Kauffman bracket skein module) 
shows that the sequences $\d^*_K(n)$ and $\d_K(n)$ 
have at most quadratic growth rate with respect to $n$. More precisely,
for {\em every} knot $K$ we have:
\begin{equation}
\lbl{eq.jbounds}
-\frac{1}{2} c^-_K n^2 +O(n) \leq \d^*_K(n) \leq  \d_K(n) \leq \frac{1}{2} c^+_K 
n^2 +O(n)
\end{equation}
This implies that the slopes $s$ of the quadratic quasi-polynomial $\d_K$
satisfy $-c^-_K \leq s \leq  c^+_K$. Replacing $K$ by its mirror $K^*$,
it implies the same inequality for the slopes $s^*$ of $\d^*_K$ and concludes
the proof of Proposition \ref{prop.c+-}.
\qed

Consider a reduced planar projection of an alternating knot $K$ with 
$c^{\pm}_K$ positive/negative crossings. Then, the number of crossings $c_K$
and the writhe $w_K$ of $K$ are given by $c_K=c^+_K+c^-_K$ and 
$w_K=c^+_K-c^-_K$. Let $\s_K$ denote the signature of $K$. Then we can
express the minimum and maximum degrees $\d^*_K(n)$ and $\d_K(n)$ of $K$
in terms of $w_K,c_K$ and $\s_K$. This
was shown by Kauffman, Murasugi and Thistlethwaite, \cite{Kf,Mu,Th}, 
and further
simplified by Turaev \cite{Tu1}. See also \cite[p.42]{Li} and 
\cite[Prop.2.1]{Le}. With our conventions, Proposition 2.1 of \cite{Le}
states that for all $n$ we have:

\begin{eqnarray}
\lbl{eq.mK}
\lbl{eq.MK}
\d_K(n) &=& \frac{(c_K+w_K)}{4} n^2 + \frac{-|A|+2c^+_K+1}{2} n 
\\
\d^*_K(n) &=& \frac{(-c_K+w_K)}{4} n^2 + \frac{|B|-2c^-_K-1}{2} n 
\end{eqnarray} 
where $|A|$ (resp. $|B|$) is the number of circles of the $A$ (resp. $B$)
smoothing of the planar projection. For example, for the right-handed 
trefoil $T$, we have
\begin{eqnarray*}
J_{T,0}(q) &=& 1 \\
J_{T,1}(q) &=& q + q^3 - q^4 \\
J_{T,2}(q) &=& q^2 + q^5 - q^7 + q^8 - q^9 - q^{10} + q^{11} \\ 
J_{T,3}(q) &=& q^3 + q^7 - 
q^{10} + q^{11} - q^{13} - q^{14} + q^{15} - q^{17} + q^{19} + q^{20} - q^{21} \\
J_{T,4}(q) &=& q^4 + 
q^9 - q^{13} + q^{14} - q^{17} - q^{18} + q^{19} - q^{22} - q^{23} + 
 2 q^{24} - q^{28} + 2 q^{29} - q^{32} - q^{33} + q^{34}
\end{eqnarray*}
and
$$
\d_T(n) = \frac{3}{2} n^2 + \frac{5}{2} n, \qquad
\d^*_T(n) = n 
$$
$$
\d^+_T(n)=\frac{3}{2} n^2 + \frac{3}{2} n, \qquad
\d^-_T(n)=\frac{3}{2} n^2 + \frac{7}{2} n
$$
and
$$
c^+_T=3, \qquad c^-_T=0, \qquad c_T=3, \qquad w_T=3, \qquad
\s_T=-2, \qquad |A|=2, \qquad |B|=3.
$$
Murasugi and Turaev observe that  \cite[p.219-220]{Tu1}

\begin{equation}
\lbl{eq.AB}
|A|+|B|=c_K+2, \qquad c_K=c^+_K + c^-_K, \qquad w_K=c^+_K-c^-_K,
\qquad \s_K=|A|-1-c^+_K=-|B|+1+c^-_K.
\end{equation}
Equation \eqref{eq.MK} implies that $\d_K(n)$ is a quadratic polynomial
(i.e., a quasi-polynomial of period $1$) with coefficient of $n^2$ equal
to $c^+_K/2$, i.e., with slope $c^+_K$. This concludes the proof of Theorem 
\ref{thm.1alt}.
Equations \eqref{eq.delta}, \eqref{eq.mK}, \eqref{eq.MK} and \eqref{eq.AB}
prove the first part of Theorem \ref{thm.2alt}.

It remains to show that the two checkerboard surfaces of a reduced projection
of an alternating knot $K$ have slopes $2 c^+_K$ and $-2 c^-_K$. Observe that
if $s=p m + q l$ is the slope of a surface $S$ (where $(m,l)$ is the
standard meridian-longitude pair) and $\la \cdot, \cdot \ra$ denotes the 
form in the boundary of a neighborhood of $K$,
then $q=\la m,s \ra$, $p=\la s,l \ra$. If $S$ is a black surface with slope 
$s=p m + q l$, then the geometrically $s$ and $m$ intersect at a point,
thus $q=\pm 1$. In addition, $s$ follows the knot $K$ as we move towards
the crossing, and intersects $l$ twice around each positive crossing,
and none around each negative crossing. The result follows.
\qed

\begin{remark}
\lbl{rem.involution}
Let $V$ denote the 3-dimensional $\BQ$-vector space spanned by the 
functions $c,w,\s$ on the set of alternating knots. There is an involution
$K \mapsto K^*$ on this set, which includes an involution on $V$: 
$$
c^*=c, \qquad w^*=-w, \qquad \s^*=-\s.
$$
On the other hand, $\d$, $\d^*$ and $\d^{\pm}$ belong to $V$ and 
$$
(\d^+)^*=\d^+, \qquad (\d^-)^*=-\d^-.
$$
Thus, $\d^+$ is a $\BQ$-linear combination of $c$, and 
$\d^-$ is a $\BQ$-linear combination of $w$ and $\s$. This is precisely the
content of Theorem \ref{thm.2alt}.
\end{remark}

\begin{remark}
\lbl{rem.degreej}
Let $f[k]$ denote the coefficient of $n^k$ in a polynomial $f(n)$.
Equations \eqref{eq.deltap} and \eqref{eq.deltan} imply that for all
alternating knots $K$ we have:
\begin{equation}
\lbl{eq.a1}
c_K = 2 \d^+_K[1]= 2 \d^+_K[2] 
\end{equation}
and
\begin{equation}
\lbl{eq.a2}
\s_K=\d^-_K[2]-\d^-_K[1], \qquad w_K=2 \d^-_K[2].
\end{equation}
\end{remark}

\section{Computing the Jones slopes and the Jones period of a knot}
\lbl{sec.nonalt}

\subsection{Some lemmas on quasi-polynomials}
\lbl{sub.quasi}

To better present the experimental (and in some cases, proven) data presented
in the next section, let
us give some lemmas on quasi-polynomials. If $a(n)$ is a sequence of numbers,
consider the generating series
\begin{equation}
\lbl{eq.Ga}
G_a(z)=\sum_{n=0}^\infty a(n) z^n.
\end{equation}
The next well-known lemma characterizes quasi-polynomials.
appears in \cite[Prop.4.4.1]{St} and \cite[Lem.3.24]{BR}.

\begin{lemma}
\lbl{lem.quasi}\cite[Prop.4.4.1]{St}\cite[Lem.3.24]{BR}
The following are equivalent:
\newline
\rm{(a)}
$a(n)$ is a quasi-polynomial of period $\pi$
\newline
\rm{(b)} the generating series 
$$
G_a(z)=\frac{P(z)}{Q(z)}
$$ 
is a rational function where $P(z), Q(z) \in \BC[z]$, every zero $\a$ of 
$Q(z)$ satisfies $\a^\pi=1$ 
(provided that $P(z)/Q(z)$ has been reduced to lowest
term) and $\mathrm{deg}P < \mathrm{deg}Q$.
\newline
\rm{(c)}
For all $n \geq 0$, 
\begin{equation}
\lbl{eq.api}
a(n)=\sum_{i=1}^k p_i(n) \ga_i^n
\end{equation}
where each $p_i(n)$ is a polynomial function of $n$ and each $\ga_i$
satisfies $\ga_i^\pi=1$.
\newline
Moreover, the degree of $p_i(n)$ in \eqref{eq.api} is one less than the
multiplicity of the root $\ga_i^{-1}$ in $Q(z)$, provided $P(z)/Q(z)$ has
been reduced to lowest terms.
\end{lemma}

\begin{definition}
\lbl{def.mono}
We say that a quadratic quasi-polynomial $a(n)$ is {\em mono-sloped} if it
has only one slope $s$. In other words, we have
$$
a(n)=\frac{s}{2} n^2 + b(n)
$$
where $b(n)$ is a linear quasi-polynomial.
\end{definition}
To phrase our next corollary, let $\Phi_n(z)$ denote
the $n$-th {\em cyclotomic polynomial}, and let $\phi(n)$ denote {\em Euler's 
$\phi$-function}.

\begin{corollary}
\lbl{cor.quasi}
$a(n)$ is mono-sloped if and only if
$$
G_a(z)=\frac{a z^2+ b z + c}{(1-z)^3} + \sum_{d>1} \frac{R_d(z)}{\Phi_d(z)^{c_d}}
$$
where the summation is over a finite set of natural numbers and $c_p \leq 2$ 
and $R_d(z) \in \BC[z]$ has degree less than $\phi(d) c_d$. Moreover,
$$
s=a+b+c.
$$
\end{corollary}

\begin{proof}
Observe that
$$
\sum_{n=0}^\infty (\a n^2 + \b n + \ga) z^n=\frac{a z^2+ b z + c}{(1-z)^3}
$$
if and only if 
$$
\a=\frac{1}{2}(a + b + c), \qquad 
\b = \frac{1}{2} (-a + b + 3 c), \qquad \ga = c.
$$
\end{proof}

\begin{proof}(of Lemma \ref{lem.spi})
Given $s$, there exists an arithmetic progression $\pi n +k$ such that
for all natural numbers $n$ we have
$$
a(\pi n+k)=\frac{s}{2} (\pi n+k)^2 + \b (\pi n+k) + \ga
$$
Now, let
$$
b(n) = \frac{s}{2} (\pi n+k)^2 + \b (\pi n+k) + \ga=
s \pi^2 \binom{n}{2} + \left(\b\pi+k \pi s+ \frac{\pi^2 s}{2}\right)n
+\ga+\b k + \frac{k^2 s}{2}.
$$
Now, $b$ takes integer values at all integers. This implies that
$$
s \pi^2, \quad 
\b\pi+k \pi s+ \frac{\pi^2 s}{2} , 
\quad
\ga+\b k + \frac{k^2 s}{2} \in  \BZ
$$
(see \cite{BR}). The result follows.
\end{proof}

It is often easier to detect the periodicity properties of the
the {\em difference} 
$$
(\De a)(n):=a(n+1)-a(n)
$$ 
of a sequence $a(n)$. It is easy to recover $G_a(z)$ from $G_{\De a}(z)$ 
and $G_a(0)$.

\begin{lemma}
\lbl{lem.DG}
With the above conventions, we have:
\begin{equation}
G_{\De a}(z)=\frac{G_a(z)(1-z)-G_a(0)}{z}
\end{equation}
\end{lemma}

\begin{proof}
We have:
\begin{eqnarray*}
G_{\De a}(z) &=& \sum_{n=0}^\infty (a(n+1)-a(n))z^n \\
&=& 1/z \sum_{n=0}^\infty a(n+1) z^{n+1} -\sum_{n=0}^\infty a(n) z^n \\
&=& 1/z(G_a(z)-G_a(0))-G_a(z).
\end{eqnarray*}
\end{proof}
We can iterate the above by considering the $k$-th difference defined
by $\De^0 a=a$ and $\De^k a=\De(\De^{k-1} a)$ for $k \geq 1$. 

\subsection{Computing the colored Jones function of a knot}
\lbl{sub.computing}

There are several ways to compute the colored Jones function $J_{K,n}(q)$
of a knot $K$. For example, one may use a planar projection and 
$R$-{\em matrices}; see for example, \cite{Tu2,GL1} and also \cite{B-N}). 
Alternatively, one may use planar projections and {\em shadow formulas}
as discussed at length in \cite{Co} and \cite{GVa}. Or one may use {\em fusion}
quantum spin networks and recoupling theory, discussed in 
\cite{CFS,KL,Co,GVa}. 
All these approaches gives various useful formulas for $J_{K,n}(q)$ presented
as a finite sum of a proper $q$-hypergeometric summand \cite{GL1}.
A careful inspection of the summand allows in several cases to compute the
degree of $J_{K,n}(q)$.

\subsection{Guessing the colored Jones function of a knot}
\lbl{sub.guess}

In this section we guess the sequence $\d_K$ of knots with a small number
of crossings, using the following strategy, inspired by conversations with
D. Zagier. 
Using the {\tt KnotAtlas} \cite{B-N} we compute as many values $J_{K,n}(q)$
of the colored Jones function as we can, and record their degree. This
gives us a table of values of the quadratic quasi-polynomials $\d_K(n)$ and 
$\d^*_K(n)$. Taking the third difference of this table results into a degree $0$
quasi-polynomial, i.e., a periodic function. At this point, we make a guess
for this periodic function, and the corresponding generating series. Then,
we use Lemma \ref{lem.DG} and our guess for the second different to obtain
a formula for $G_{\d_K}(z)$ and $G_{\d^*_K}(z)$. The partial fraction decomposition
then gives us a formula for $\d_K(n)$ and  $\d^*_K(n)$. 
In some cases, using explicit finite multi-dimensional sum formulas for the 
colored Jones 
polynomial, one can prove that the guessed formula for $\d_K(n)$ and  
$\d^*_K(n)$ are indeed correct. In this section, we will not bother with
proofs.

As an example of our method, we will guess a formula for $\d(n)$ and $\d^*(n)$
for the $(-2,3,7)$ pretzel knot. One can actually prove that our guess
is correct, using the fusion formulas for the 3-pretzel knots, but we will
not bother. The values of $\d^*(n)$ starting with $n=0$ are given by:
\begin{eqnarray*}
\d^* & : & 0, 5, 10, 15, 20, 25, 30, 35, 40, 45, 50, 55, 60, 65, 70, \dots
\end{eqnarray*}
Taking the first difference we get the following values of $(\D \d^*)(n)$
starting with $n=0$:
\begin{eqnarray*}
\D \d^* & : & 5,5,5,5,5,5,5,5,5,5,5,5,5,5,\dots
\end{eqnarray*}
appears to be the constant sequence from which
we guess that $\d^*(n)=5n$, and correspondingly the generating series
is
$$
G_{\d^*}(z)=\frac{5 z}{(1-z)^2}
$$
More interesting is the sequence $\d(n)$ starting with $n=0$:
\begin{eqnarray*}
\d & : & 0, 13, 35, 67, 108, 158, 217, 286, 364, 451, 547, 653, 768, 892, 
1025, 1168, 1320, 1481, 1651, 1831, \dots
\end{eqnarray*}
Taking the first, second and third difference we obtain:
\begin{eqnarray*}
\De \d 
& : & 13, 22, 32, 41, 50, 59, 69, 78, 87, 96, 106, 115, 124, 133, 143, 152, \
161, 170, 180, \dots 
\\
\De^2 \d 
& : & 9, 10, 9, 9, 9, 10, 9, 9, 9, 10, 9, 9, 9, 10, 9, 9, 9, 10, \dots
\\
\De^3 \d 
& : & 1, -1, 0, 0, 1, -1, 0, 0, 1, -1, 0, 0, 1, -1, 0, 0, 1, \dots
\end{eqnarray*}
Thus, we guess that $\De^3 \d$ is a periodic sequence with period $4$ and 
generating series
$$
G_{\De^3 \d}(z)=\sum_{n=0}^\infty \left(z^{4n}-z^{4n+1}\right)=\frac{1}{(1+z)(1+z^2)}
$$
Using Lemma \ref{lem.DG} three times, we compute
$$
G_{\d}(z)=\frac{13 z + 9 z^2 + 10 z^3 + 9 z^4 - 4 z^5}{
(1 - z)^3 (1 + z + z^2 + z^3)}=
\frac{-3 + 216 z - 65 z^2}{16 (1 - z)^3} +
\frac{3 + 4 z - z^2}{16 (1 + z + z^2 + z^3)}
$$
Taking the partial fraction decomposition, it follows that
$$
\d(n)=\left[\frac{37}{8} n^2 + \frac{17}{2} n \right]=
\frac{37}{8} n^2 + \frac{17}{2} n + \e(n),
$$
where $\e(n)$ is a periodic sequence of period $4$ given by:
$$
\e(n)=\begin{cases}
0 & \text{if} \qquad n=0 \bmod 4 \\
\frac{1}{8} & \text{if} \qquad n=1 \bmod 4 \\
\frac{1}{2} & \text{if} \qquad n=2 \bmod 4 \\
\frac{1}{8} & \text{if} \qquad n=3 \bmod 4 \\
\end{cases}
$$

\subsection{A summary of non-alternating knots}
\lbl{sub.summary}

In this section we list the quasi-polynomials $\d_K$ and $\d^*_K$ 
of non-alternating knots $K$ with $8$ and $9$ crossings.
In the Rolfsen table of knots, the non-alternating knots with
$8$ crossings are $8_k$ where 
$k=19,\dots,21$, with $9$ crossings are $9_k$ where $k=42,\dots,49$.
Let us give a combined table of the non-alternating knots $K$ 
with $8$ and $9$ crossings, their period $\pi_K$, their Jones
slopes $\js_K$ and $\js^*_K$, and their distinct boundary slopes.
The boundary slopes $\bs_K$ are computed using the program of \cite{HO},
corrected in \cite{Du}, which computes the boundary slopes of 
all Montesinos knots except $9_{49}$ which is not a Montesinos knots.
In all those cases, the set of boundary slopes agrees with the slopes
of the $A$-polynomial of \cite{Cu,CCGLS}, once $0$ is included.

\begin{center}
\begin{tabular}{||c|c|c|c|c||}
\hline
$ K    $ &  $ \pi_K $ & $  \js_K $ & $ \js^*_K     $ & $ \bs_K     $ \\
\hline
$ 8_{19} $ & $ 2  $ & $ 6 $ & $ 0       $ & $   \{0,12\} $ \\
$ 8_{20} $ & $ 3  $ & $ 4/3 $ & $  -5   $ & $ \{-10, 0, 8/3\}    $ \\
$ 8_{21} $ & $ 2  $ & $ 1/2  $ & $  -6  $ & $   \{-12, -6, -2, 0, 1\}$ \\
\hline
$ 9_{42} $ & $ 2  $ & $ 3 $ & $ -4      $ & $   \{-8, 0, 8/3, 6\}$ \\
$ 9_{43} $ & $ 3  $ & $ 16/3 $ & $ -2   $ & $ \{-4, 0, 6, 8, 32/3\}$ \\
$ 9_{44} $ & $ 3  $ & $ 7/3 $ & $ -5    $ & $   \{-10, -2, 0, 1, 2, 14/3\}$ \\
$ 9_{45} $ & $ 2  $ & $ 1/2 $ & $  -7   $ & $ \{-14, -10, -8, -4, -2, 0, 1\}$ \\
$ 9_{46} $ & $ 2  $ & $ 1 $ & $  -6     $ & $   \{-12, 0, 2\} $ \\
$ 9_{47} $ & $ 2  $ & $ 9/2 $ & $  -3   $ & $   \{-6, 0, 4, 8, 9, 16\} $ \\
$ 9_{48} $ & $ 2  $ & $ 11/2 $ & $  -2  $ & $   \{-4, 0, 4, 8, 11\} $ \\
$ 9_{49} $ & $ 2  $ & $ 15/2 $ & $ 0    $ & $   \{0, 4, 6, 12, 15\} $ \\ 
\hline
\end{tabular}
\end{center}

The above data are in agreement with Conjecture \ref{conj.1}.
Let us make a phenomenological remark regarding all examples 
of non-alternating knots with $8$ or $9$ crossings.

\begin{itemize}
\item[(a)]
$\d^*(n)$ and $\d(n)$ are mono-sloped, i.e., they are of the form 
$s n^2/2+\e(n)$
where $\e(n)$ is a linear quasi-polynomial.
\item[(b)]
For all knots, $2 \js$ is a boundary slope, though not necessarily the 
largest one.
\item[(c)]
In the case of the $8_{20}, 9_{43}$ and $9_{44}$ knots, the degree of 
$\e(n)$ is $1$ and in all other cases, it is zero.
\item[(d)]
The period of all non-alternating knots 
was greater than $1$. $8_{20}, 9_{43}, 9_{44}$ knots have period $3$, and
$(-2,3,7)$ has period $4$. The period of $(-2,3,p)$ for odd $p \geq 5$
appears  to be $p-3$, and the number of crossings is $p+5$. Thus the period
can be asymptotically as large as the number of crossings. 
\item[(e)]
For the case of $8_{21},9_{45},9_{46},9_{47}$, the leading coefficient 
is $2$ and for $9_{48},9_{49}$ it is $-2$. 
\end{itemize}

\subsection{The $8$-crossing non-alternating knots}
\lbl{sub.8na}

\subsubsection{The $8_{19}$ knot}
\lbl{sub.819}

In the data below, we will give the first few values of $\d^*(n)$ and $\d(n)$,
the guessed decomposition of the generating series $G_{\d^*}(z)$ and $G_{\d}(z)$
of the quasi-polynomials $\d^*$ and $\d$.

Some values of $\d^*(n)$ and $\d(n)$ starting with $n=0$:
\begin{eqnarray*}
\d^* & : & 0, 3, 6, 9, 12, 15, 18, 21, \dots
\\
\d & : & 
0, 8, 23, 43, 70, 102, 141, 185, \dots
\end{eqnarray*}
\begin{eqnarray*}
G_{\d^*}(z) & = & \frac{3z}{(1-z)^2}
\\
G_{\d}(z) & = & \frac{8 z + 7 z^2 - 3 z^3}{(1 - z)^3 (1 + z)} 
\\ 
& = &
\frac{-1 + 36 z - 11 z^2}{4 (1 - z)^3} + \frac{1}{4 (1 + z)}
\end{eqnarray*}
\begin{eqnarray*}
\d^*(n) & = & 3n
\\
\d(n) & = & 3 n^2 + \frac{11}{2} n -\frac{1}{4}+ \e(n)
\end{eqnarray*}
where $\e(n)=(-1)^n/4$ is a 2-periodic sequence.
Note that in this example the values of $\d^*_K(n)$ for $n=0,1,2,3$ suffice to
prove that $\d_K(n)$ is {\em not} a polynomial of $n$.

\subsubsection{The $8_{20}$ knot}
\lbl{sub.820}

\begin{eqnarray*}
\d^* & : & 0, -5, -15, -30, -50, -75, -105, -140, -180, -225, -275, -330, -390, 
-455, -525, -600, -680, 
\\
& & -765, -855, -950, -1050, \dots
\\
\d & : & 0, 1, 2, 7, 12, 16, 26, 35, 42, 57, 70, 80, 100, 117, 130, 155, 176, 
192, 222, 247, 266 , \dots
\end{eqnarray*}
\begin{eqnarray*}
G_{\d^*}(z) & = & -\frac{5 z}{(1 - z)^3}
\\
G_{\d}(z) & = & \frac{z + z^2 + 5 z^3 + 3 z^4 + 2 z^5}{(1 - z)^3 (1 + z + z^2)^2}
\\ 
& = &
\frac{-2 + 12 z + 2 z^2}{9 (1 - z)^3}
+\frac{2 + 7 z + 4 z^2 + 2 z^3}{9 (1 + z + z^2)^2}
\end{eqnarray*}
\begin{eqnarray*}
\d^*(n) & = & -\frac{5 n(n+1)}{2}
\\
\d(n) & = & \frac{2}{3} n^2 + \frac{2}{9} n -\frac{2}{9}+\ep(n)
\end{eqnarray*}
where $\ep(n)$ is a {\em linear} quasi-polynomial with period $3$.

\subsubsection{The $8_{21}$ knot}
\lbl{sub.821}

\begin{eqnarray*}
\d^* & : & 0, -7, -20, -39, -64, -95, -132, -175, \dots
\\
\d & : & 0, -1, -1, -1, 0, 1, 3, 5 , \dots
\end{eqnarray*}
\begin{eqnarray*}
G_{\d^*}(z) & = & \frac{-7 z + z^2}{(1 - z)^3}
\\
G_{\d}(z) & = & \frac{-z + z^2 + z^3}{(1 - z)^3 (1 + z)}
\\ 
& = &
\frac{-1 - 4 z + 9 z^2}{8 (1 - z)^3}+\frac{1}{8 (1 + z)}
\end{eqnarray*}
\begin{eqnarray*}
\d^*(n) & = & -n(3n+4)
\\
\d(n) & = & \frac{1}{4} n^2 - n -\frac{1}{8}-\e(n)
\end{eqnarray*}
where $\e(n)=(-1)^n/8$ is a 2-periodic sequence.

\subsection{The $9$-crossing non-alternating knots}
\lbl{sub.9na}

\subsubsection{The $9_{42}$ knot}
\lbl{sub.942}

\begin{eqnarray*}
\d^* & : & 0, -3, -10, -21, -36, -55, -78, -105, \dots
\\
\d & : & 0, 3, 10, 19, 32, 47, 66, 87, \dots
\end{eqnarray*}
\begin{eqnarray*}
G_{\d^*}(z) & = & \frac{-3 z - z^2}{(1 - z)^3}
\\
G_{\d}(z) & = & \frac{z (3 + 4 z - z^2)}{(1 - z)^3 (1 + z)}
\\ 
& = &
\frac{1 - 16 z + 3 z^2}{4 (-1 + z)^3} + \frac{1}{4 (1 + z)}
\end{eqnarray*}
\begin{eqnarray*}
\d^*(n) & = & -n (2 n + 1)
\\
\d(n) & = & \frac{3}{2} n^2 + 2 n -\frac{1}{4} +\e(n)
\end{eqnarray*}
where $\e(n)=(-1)^n/4$ is a 2-periodic sequence.

\subsubsection{The $9_{43}$ knot}
\lbl{sub.943}

\begin{eqnarray*}
\d^* & : & 0, 0, -2, -6, -12, -20, -30, -42
\\
\d & : & 0, 7, 17, 37, 60, 85, 122, 161 
\end{eqnarray*}
\begin{eqnarray*}
G_{\d^*}(z) & = & (2 z)/(-1 + z)^3
\\
G_{\d}(z) & = & 
(z (-7 - 10 z - 20 z^2 - 9 z^3 - 5 z^4 + 3 z^5))/((-1 + z)^3 (1 + z + z^2)^2)
\\ 
& = &
(5 - 72 z + 19 z^2)/(9 (-1 + z)^3) +
(5 + 16 z + 13 z^2 + 8 z^3)/(9 (1 + z + z^2)^2)
\end{eqnarray*}
\begin{eqnarray*}
\d^*(n) & = & -n ( n - 1)
\\
\d(n) & = & -(5/9) + (38 n)/9 + (8 n^2)/3 + \e(n)
\end{eqnarray*}
where $\e(n)$ is a linear quasi-polynomial with period $3$.

\subsubsection{The $9_{44}$ knot}
\lbl{sub.944}

\begin{eqnarray*}
\d^* & : & 0, -5, -15, -30, -50, -75, -105, -140, \dots
\\
\d & : & 0, 2, 5, 13, 22, 31, 47, 63, \dots
\end{eqnarray*}
\begin{eqnarray*}
G_{\d^*}(z) & = & (5 z)/(-1 + z)^3
\\
G_{\d}(z) & = & 
-((z (2 + 3 z + 8 z^2 + 5 z^3 + 3 z^4))/((-1 + z)^3 (1 + z + z^2)^2))
\\ 
& = &
(2 - 21 z - 2 z^2)/(9 (-1 + z)^3) +
(2 + 7 z + 4 z^2 + 2 z^3)/(9 (1 + z + z^2)^2) 
\end{eqnarray*}
\begin{eqnarray*}
\d^*(n) & = & -5 n ( n + 1)/2
\\
\d(n) & = & -(2/9) + (13 n)/18 + (7 n^2)/6 + \e(n)
\end{eqnarray*}
where $\e(n)$ is a linear quasi-polynomial with period $3$.

\subsubsection{The $9_{45}$ knot}
\lbl{sub.945}

\begin{eqnarray*}
\d^* & : & 0, -8, -23, -45, -74, -110, -153, -203, \dots
\\
\d & : & 0, -1, -1, -1, 0, 1, 3, 5, \dots
\end{eqnarray*}
\begin{eqnarray*}
G_{\d^*}(z) & = & -((-8 z + z^2)/(-1 + z)^3)
\\
G_{\d}(z) & = & -((z (-1 + z + z^2))/((-1 + z)^3 (1 + z)))
\\ 
& = &
(1 + 4 z - 9 z^2)/(8 (-1 + z)^3) +
1/(8 (1 + z))
\end{eqnarray*}
\begin{eqnarray*}
\d^*(n) & = & -n (7 n + 9)/2
\\
\d(n) & = & -(1/8) - n + n^2/4 + \e(n)
\end{eqnarray*}
where $\e(n)=(-1)^n/8$ is a periodic sequence with period $2$.

\subsubsection{The $9_{46}$ knot}
\lbl{sub.946}

\begin{eqnarray*}
\d^* & : & 0, -6, -18, -36, -60, -90, -126, -168, \dots
\\
\d & : & 0, 0, 2, 4, 8, 12, 18, 24, \dots
\end{eqnarray*}
\begin{eqnarray*}
G_{\d^*}(z) & = & -\frac{6 z}{(1 - z)^3}
\\
G_{\d}(z) & = & \frac{2 z^2}{(1 - z)^3 (1 + z)}
\\ 
& = &
\frac{-1 + 4 z + z^2}{4 (1 - z)^3}
+ \frac{1}{4 (1 + z)}
\end{eqnarray*}
\begin{eqnarray*}
\d^*(n) & = & -3n(n+1)
\\
\d(n) & = & -(1/4) + n^2/2 + \e(n)
\end{eqnarray*}
where $\e(n)=(-1)^n/4$ is a periodic sequence with period $2$.

\subsubsection{The $9_{47}$ knot}
\lbl{sub.947}

\begin{eqnarray*}
\d^* & : & 0, -2, -7, -15, -26, -40, -57, -77, \dots
\\
\d & : & 0, 5, 15, 29, 48, 71, 99, 131, \dots
\end{eqnarray*}
\begin{eqnarray*}
G_{\d^*}(z) & = & (2 z + z^2)/(-1 + z)^3
\\
G_{\d}(z) & = & (z (-5 - 5 z + z^2))/((-1 + z)^3 (1 + z))
\\ 
& = &
(1 - 44 z + 7 z^2)/(8 (-1 + z)^3) +
1/(8 (1 + z))
\end{eqnarray*}
\begin{eqnarray*}
\d^*(n) & = & -n (3 n + 1)/2
\\
\d(n) & = &  -(1/8) + 3 n + (9 n^2)/4 + \e(n)
\end{eqnarray*}
where $\e(n)=(-1)^n/8$ is a periodic sequence with period $2$.

\subsubsection{The $9_{48}$ knot}
\lbl{sub.948}

\begin{eqnarray*}
\d^* & : & 0, -1, -4, -9, -16, -25, \dots
\\
\d & : & 0, 6, 18, 35, 58, 86, \dots
\end{eqnarray*}
\begin{eqnarray*}
G_{\d^*}(z) & = & -((-z - z^2)/(-1 + z)^3)
\\
G_{\d}(z) & = & (z (-6 - 6 z + z^2))/((-1 + z)^3 (1 + z))
\\ 
& = &
(1 - 52 z + 7 z^2)/(8 (-1 + z)^3) + 
1/(8 (1 + z))
\end{eqnarray*}
\begin{eqnarray*}
\d^*(n) & = & -n^2
\\
\d(n) & = & -(1/8) + (7 n)/2 + (11 n^2)/4 + \e(n)
\end{eqnarray*}
where $\e(n)=(-1)^n/8$ is a periodic sequence with period $2$.

\subsubsection{The $9_{49}$ knot}
\lbl{sub.949}

\begin{eqnarray*}
\d^* & : & 0, 2, 4, 6, 8, 10, \dots
\\
\d & : & 0, 9, 26, 50, 82, 121, \dots
\end{eqnarray*}
\begin{eqnarray*}
G_{\d^*}(z) & = & (2 z)/(-1 + z)^2
\\
G_{\d}(z) & = & (z (-9 - 8 z + 2 z^2))/((-1 + z)^3 (1 + z))
\\ 
& = &
(1 - 76 z + 15 z^2)/(8 (-1 + z)^3) +
1/(8 (1 + z))
\end{eqnarray*}
\begin{eqnarray*}
\d^*(n) & = & 2 n
\\
\d(n) & = & -(1/8) + (11 n)/2 + (15 n^2)/4 + \e(n)
\end{eqnarray*}
where $\e(n)=(-1)^n/8$ is a periodic sequence with period $2$.

\subsection{The case of the $(-2,3,p)$ pretzel knots}
\lbl{sub.pretzel}

A triple sum formula for the colored Jones polynomial of pretzel knots
with 3 pretzels is available, and using it we can compute $\d_K(n)$
and $\d^*_K(n)$ for all pretzel knots of the form $(-2,3,p)$ for odd $p$;
see \cite{GVa}. 
We will state the result of the computation here. 
Recall the $k$-th difference $\De^k f$ of a sequence $f$
from Section \ref{sub.quasi}. When $p>0$ is odd, we have:

$$
G_{\d^*}(z) = \frac{(p+3)z}{2(1-z)^2},
\qquad
G_{\De^3 \d}(z) =
\begin{cases} 
\frac{z^{p-7}(1-z)}{1-z^{p-3}} & \qquad p \geq 7 \\
 -\frac{3}{1+z} & \qquad p=5 \\
 -\frac{2}{1+z} & \qquad p=3 \\
 0 & \qquad p=1 
\end{cases}
$$
It follows that
$$
\d^*(n)= \frac{(p+3)n}{2}, \qquad 
\d(n)=\begin{cases}
\frac{5}{2} n^2 + \e_p(n)  & \qquad p=1 \\
3 n^2 + \e_p(n)  & \qquad p=3 \\
\frac{(p^2-p-5)n^2}{2(p-3)} + \e_p(n) & \qquad p \geq 5
\end{cases}
$$
where $\e_p(n)$ are linear quasi-polynomials. When $p<0$ is odd we have:
$$
\De^3 G_{\d^*}(z) =
\begin{cases}
0 & \qquad p=-1 \\
-\frac{4+4z+3z^2+z^3}{(1+z+z^2)^2} & \qquad p=-3 \\
\frac{z^{|p|-4}-2 z^{|p|-3} -\sum_{k=|p|-2}^{2|p|-4} z^k}{(\sum_{k=0}^{|p|-1}
z^k)^2} & \qquad p \leq -5
\end{cases}
\qquad
G_{\d}(z) = \frac{z(p+13-(p+3)z)}{2(1-z)^3}
$$
It follows that 
$$
\d^*(n)=
\begin{cases}
\frac{5}{2} n^2 + \e_p(n)   & \qquad p = -1 \\
\frac{(p+1)^2 n^2}{2 p}+\e_{p}(n) & \qquad p \leq -3
\end{cases}
\qquad
\d(n)=\frac{n(5n+(p+8))}{2}
$$
where $\e_p(n)$ are linear quasi-polynomials.
Notice that the above formulas single out exceptional behavior at
$p=-3,-1,1,3,5$.
The Jones period and the Jones slopes are given by 
\begin{equation}
\lbl{eq.jspretzel}
\pi=\begin{cases}
p-3 & \quad p \geq 5 \\
2 & \quad p=3 \\
|p| & \quad p \leq 1 
\end{cases}
\qquad \js=
\begin{cases}
\frac{p^2-p-5}{p-3} & \quad p \geq 5 \\
6 & \quad p=3 \\
5 & \quad p \leq -1 
\end{cases}
\qquad \js^*=
\begin{cases}
0 & \quad p \geq 5 \\
0 & \quad p=3 \\
\frac{(p+1)^2}{p} & \quad p \leq 1 
\end{cases}
\end{equation}
On the other hand, Hatcher-Oertel and Dunfield (see \cite{HO} and \cite{Du}) 
compute the slopes of those Montesinos knots
\begin{equation}
\lbl{eq.bspretzel}
\bs_p=
\begin{cases}
\{0, 16, \frac{2 (p^2 -p -5)}{p -3}, 2 (3 + p)\} & p \geq 7 \\
\{0,10, 2\frac{(p+1)^2}{p}, 2(p+3) \} & p \leq -1
\end{cases}
\end{equation}
Compare also with Mattman \cite[p.32]{Ma} who computes which of those
slopes are visible from the geometric component of the $A$-polynomial.
Equations \eqref{eq.jspretzel} and \eqref{eq.bspretzel} together with
the fact that $0$ is a boundary slope confirm
Conjecture \ref{conj.1} for all $(-2,3,p)$ pretzel knots.

\subsection{The case of torus knots}
\lbl{sub.torus}

In this section we will use Morton's formula for the colored Jones function
of a torus knot to compute the degree of the colored Jones function
and verify Conjecture \ref{conj.1}. 

Let $T(a,b)$ denote the $(a,b)$ torus knot for a pair of coprime integers
$a,b$. Since the mirror image of $T(a,b)$ is $T(a,-b)$, we will focus
on the case of $a,b >0$. With our conventions, Morton's formula \cite{Mo}
for the colored Jones function is the following:

\begin{equation}
\lbl{eq.morton}
J_{T(a,b),n}(q)=\frac{q^{\frac{1}{4}ab n(n+2)}}{q^{\frac{n+1}{2}}-q^{-\frac{n+1}{2}}}
\sum_{k=-\frac{n}{2}}^{\frac{n}{2}} \left(q^{-abk^2+(a-b)k+\frac{1}{2}}-q^{-abk^2+(a+b)k-\frac{1}{2}}
\right)
\end{equation}
For example, 
\begin{eqnarray*}
J_{T(2,3),1}(q)&=& q + q^3 - q^4 \\
J_{T(3,4),2}(q)&=& q^6 + q^9 + q^{12} - q^{13} - q^{16} - q^{19} + q^{20} - q^{22} + q^{23}
\end{eqnarray*}
The summand of Equation \eqref{eq.morton} consists of two monomials with
exponents quadratic functions of $k$. A little calculation reveals that
the maximum and minimum degree of the colored Jones function is given by
\begin{eqnarray*}
\d_{T(a,b)}(n)&=& \frac{a b}{4} n^2 + \frac{a b - 1}{2} n - (1 - (-1)^n)  
\frac{(a - 2) (b - 2)}{8} \\
\d^*_{T(a,b)}(n)&=& \frac{(a - 1) (b - 1)}{2} n 
\end{eqnarray*}
Thus the period $\pi_{T(a,b)}$ is $2$ when $a, b \neq 2$ and $1$ when
$a=2$ or $b=2$. Since $T(a,b)$ is alternating if and only if $a=2$ or $b=2$,
our results are in accordance to Question \ref{que.1}.
The boundary slopes of $T(a,b)$ are $\{0,ab\}$ (see
for example \cite{HO}). This confirms Conjecture \ref{conj.1} for torus
knots.

\subsection{Acknowledgment}
The author wishes to thank K. Ichihara, Y. Kabaya, T. Mattman,
P. Ozsv\'ath, A. Shumakovitch, R. van der Veen, D. Zagier, D. Zeilberger
and especially 
M. Culler and N. Dunfield for numerous enlightening conversations.
The idea of the Slope Conjecture was conceived during the New York
Conference on {\em Interactions between Hyperbolic Geometry, 
Quantum Topology and Number Theory} in New York in the summer of 2009.
The author wishes to thank the organizers, A. Champanerkar, O. Dasbach, 
E. Kalfagianni, I. Kofman, W. Neumann and N. Stoltzfus for their
hospitality.




\ifx\undefined\bysame
        \newcommand{\bysame}{\leavevmode\hbox
to3em{\hrulefill}\,}
\fi

\end{document}